# Moebius functions
# Irreducible Polynomials and
# Dirichlet Series

*Dang Vu Giang*
**Hanoi Institute of Mathematics**
**18 Hoang Quoc Viet**
**10307 Hanoi Vietnam**
dangvugiang@yahoo.com

***Classical Moebius function*** $\mu$ is defined on the positive integers $\mathbb{Z}_+$ and taking $0, \pm 1$ as its values. More exactly, $\mu(1) = 1$, $\mu(n) = 0$ if $n$ is divisible by the square of a prime and finally $\mu(n) = (-1)^k$ if $n$ is the product of $k$ distinct primes. Clearly, for every co-prime integers $m,n$ we have $\mu(mn) = \mu(n)\mu(m)$. We have the following beautiful inversion formulas:

**Additive formula:** Let $g : \mathbb{Z}_+ \to \mathbb{Z}$ be a function and $f(n) = \sum_{d|n} g(d)$. Then $g(n) = \sum_{d|n} \mu\left(\frac{n}{d}\right) f(d)$.

**Multiplicative formula:** Let $G$ be a commutative group, $g : \mathbb{Z}_+ \to G$ be a function and $f(n) = \prod_{d|n} g(d)$. Then $g(n) = \prod_{d|n} f(d)^{\mu\left(\frac{n}{d}\right)}$.

They can easily proved via the trivial identity $0 = \sum_{d|n} \mu(d)$ for every $n > 1$. The multiplicative formula has an interesting application for Gauss ***primitive (cyclotomic) polynomials*** $\Phi_n$. By definition $\Phi_n(x) = \prod_{(k,n)=1} (x - \omega^k)$ where $\omega$ is a $n$th primitive root of unity. For example, $\Phi_1(x) = x - 1$, $\Phi_2(x) = x + 1$, $\Phi_3(x) = x^2 + x + 1$, $\Phi_4(x) = x^2 + 1$, $\Phi_5(x) = x^4 + x^3 + x^2 + x + 1$, $\Phi_6(x) = x^2 - x + 1$. Let $\phi(n) = \sum_{(k,n)=1} 1 = n \prod_{p|n} \left(1 - \frac{1}{p}\right)$ be the Euler function. Then we have $n = \sum_{d|n} \phi(d)$. This implies that $x^n - 1 = \prod_{d|n} \Phi_d(x)$. Thus, in virtue of the multiplicative formula $\Phi_n(x) = \prod_{d|n} (x^d - 1)^{\mu(n/d)}$. Therefore, $\Phi_n$ is of integer coefficients. Moreover, $\Phi_n$ is irreducible over rationals $\mathbb{Q}$. If otherwise, there is a monic irreducible polynomial $u \in \mathbb{Z}[x]$ such that $\Phi_n = uv$, where $v \in \mathbb{Z}[x]$. Let $p$ be a prime not dividing $n$ and let $\omega$ be a root of $u$. If $\omega^p$ is not a root of $u$ then $\omega^p$ is a root of $v$. Thus, $\omega$ is root of $v(x^p)$. Consequently, $u(x) \in \mathbb{Z}[x]$ is dividing $v(x^p)$.

Reduced modulo $p$, $v(x^p) = v(x)^p$ in the polynomial ring $\mathbb{F}_p[x]$. Consequently, $u(x) \in \mathbb{F}_p[x]$ is dividing $v(x) \in \mathbb{F}_p[x]$. Therefore, in an extension of $\mathbb{F}_p$ the polynomial $x^n - 1$ has a multiplicative root. This is a contradiction, because $p$ is not a divisor of $n$. Thus, $\omega^p$ is a root of $u$. Hence, $\omega^k$ is root of $u$ for any integer $k$ relatively prime to $n$. Therefore, $u = \Phi_n$ so $\Phi_n$ is irreducible over rationals $\mathbb{Q}$. This is known as minimal polynomials of any primitive $n$th root of 1. Now let $I_q(n)$ be the **number of monic irreducible polynomials** $f$ of degree $n$ over the **finite field** $\mathbb{F}(q) = \mathbb{F}_q$ (of $q = p^m$ elements). For example, $x^2 + x + 1$ is the only irreducible polynomial of degree 2 over $\mathbb{F}_2$. It is well known that there is $\omega \in \mathbb{F}(q)$ such that $1, \omega, \omega^2, \cdots, \omega^{m-1}$ constitute a basis for $\mathbb{F}(q)$ as vector space over $\mathbb{F}_p$. Let $\omega^m = a_1 \omega^{m-1} + \cdots + a_{m-1} \omega + a_m$ where $a_1, \cdots, a_{m-1}, a_m \in \mathbb{F}_p$. Clearly, the polynomial $h(x) = x^m - a_1 x^{m-1} - \cdots - a_{m-1} x - a_m$ is irreducible over $\mathbb{F}_p$ and and $\mathbb{F}(q) \cong \mathbb{F}_p[x]/(h)$. Hence, **two finite fields of the same cardinality are isomorphic**. Moreover, the multiplicative group $\mathbb{F}_q^*$ (of nonzero elements of $\mathbb{F}_q$) is cyclic. In fact if $\alpha \in \mathbb{F}_q^*$ is an element of maximal order $r \leq q-1$ then any element $\beta$ of $\mathbb{F}(q)^*$ satisfying $\beta^r = 1$. Indeed, if $\ell$ is the order of $\beta$ and $\pi$ is a prime then we can write $r = \pi^a r_0$ and $\ell = \pi^b \ell_0$ where $r_0$ and $\ell_0$ are not divisible by $\pi$. Clearly, $\alpha^{\pi^a}$ has order $r_0$ and $\beta^{\ell_0}$ has order $\pi^b$. Hence, $\alpha^{\pi^a} \beta^{\ell_0}$ has order $\pi^b r_0 \leq r = \pi^a r_0$. Consequently, $b \leq a$ and every divisor of $\ell$ is also a divisor of $r$. Therefore, $\ell$ itself is a divisor of $r$. Hence, $\beta^r = 1$ and the polynomial $\prod_{\vartheta \in \mathbb{F}(q)^*} (x - \vartheta)$ of degree $q - 1$ is a divisor of $x^r - 1$. Consequently, $\ell \geq q - 1$. But $\ell \leq q - 1$ so we have $\ell = q - 1$ and the multiplicative group of any finite field is cyclic. Thus, the extensions $\mathbb{F}(q)[x]/(f)$ of $\mathbb{F}(q)$ ($f$ is any irreducible polynomial over $\mathbb{F}(q)$ of degree $n$) are isomorphic and irreducible polynomials over $\mathbb{F}(q)$ of degree $n$ are completely reducible in $\mathbb{F}(q^n)$, and they are factors of $x^{q^n} - x$. More exactly, $x^{q^n} - x$ is the product of all monic irreducible polynomials of degree $d$ where $d$ is running over divisors of $n$. If otherwise, there is an irreducible factor $v$ of $x^{q^n} - x$ with degree $k$ which is not diving $n$. Clearly, $k > 1$ and $v \mid x^{q^k} - x$. But $n = k\ell + r$ with $r \in (0, k)$ so $v$ is an irreducible factor of and $x^{q^r} - x$. Thus, any element $\omega$ in the quotient field $\mathbb{F}(q)[x]/(v) = \mathbb{F}(q^k)$ is satisfying $\omega^{q^r} = \omega$ so $q^r - 1$ is divisible by $q^k - 1$. This is contradiction, because $r \in (0, k)$ and we have $x^{q^n} - x$ is the product of all monic irreducible polynomials of degree $d$ where $d$ running over divisors of $n$. Therefore, $q^n = \sum_{d \mid n} d I_q(d)$ so

$$I_q(n) = \frac{1}{n} \sum_{d|n} \mu(n/d) q^d.$$

***The famous Riemann zeta function*** $\varsigma$ is defined by $\varsigma(z) = \sum_{n=1}^{\infty} \frac{1}{n^z}$ for $\operatorname{Re}(z) > 1$. We have a reciprocal formula $\frac{1}{\varsigma(z)} = \sum_{n=1}^{\infty} \frac{\mu(n)}{n^z}$. On the other hand, $|\varsigma(z)| \leq \varsigma(\operatorname{Re} z) = \sum_{n=1}^{\infty} \frac{1}{n^{\operatorname{Re} z}} < 1 + \int_{1}^{\infty} \frac{dt}{t^{\operatorname{Re} z}} = \frac{\operatorname{Re} z}{\operatorname{Re} z - 1}$ and also $\frac{1}{|\varsigma(z)|} \leq \varsigma(\operatorname{Re} z)$ in the virtue of the reciprocal formula. Thus, $\frac{\operatorname{Re} z - 1}{\operatorname{Re} z} < |\varsigma(z)| < \frac{\operatorname{Re} z}{\operatorname{Re} z - 1}$ and consequently, $\lim_{\operatorname{Re} z \to \infty} |\varsigma(z)| = 1$ and $|\varsigma(z)| > 0$ for $\operatorname{Re}(z) > 1$. If Riemann hypothesis is not true then for every $k$ there exist a complex number $z$ with $\operatorname{Re}(z) > 1$ and an integer $N \geq k$ such that $\sum_{n=1}^{N} \frac{1}{n^z} = 0$. Now we wish to count the number of ***derangements*** of a finite set $S$, that is the number of permutations of $S$ which have no fixed points. Let $T$ be a subset of $S$. We define

$f(T) = $ the number of permutations of $S$ which fix all the elements of $T$ but fix no other element of the complement $T'$ of $T$ in $S$;

$g(T) = $ the number of permutations of $S$ which fix all the elements of $T$ but perhaps some additional elements as well. Then $g(T) = \sum_{U \supseteq T} f(U)$ where of course, $U$ is a subset of $S$. We have $g(U) = |S \setminus U|!$ and $f(\varnothing)$ is the number of derangements of $S$. Here $|P|$ denotes the cardinality of $P$. The object is to invert the formula $g(T) = \sum_{U \supseteq T} f(U)$ to obtain $f(T)$ in term of $g(U)$. To this end we define the ***Moebius function on a locally finite partially ordered set*** and will have $f(\varnothing) = |S|! \sum_{k=0}^{|S|} \frac{(-1)^k}{k!} = \left[ \frac{|S|!}{e} \right]$ where $[x]$ denote the integer closest to $x$. A ***partial order*** is a binary relation "≤" over a set $P$ which is reflexive, antisymmetric, and transitive, i.e., for all $a$, $b$, and $c$ in $P$, we have that

- $a \leq a$ (reflexivity);
- if $a \leq b$ and $b \leq a$ then $a = b$ (antisymmetry);
- if $a \leq b$ and $b \leq c$ then $a \leq c$ (transitivity).

Examples:

- Positive integers ordered by divisibility;
- Finite subsets of some set $E$, ordered by inclusion;
- Subsets of some finite set $S$, ordered by exclusion: $U \leq T$ if $U \supseteq T$.

A set with a partial order is called a partially ordered set (also called a poset). A **locally finite** poset is one for which every closed interval $[a, b] = \{x : a \leq x \leq b\}$ within it is finite. In theoretical physics a locally finite poset is also called a causal set and has been used as a model for spacetime. **Moebius function on a locally finite poset** $P$ is defined as follows: $\mu(x,x) = 1$ for all $x \in P$ and $\mu(x,y) = -\sum_{x \leq t < y} \mu(x,t)$. The local finiteness of $P$ assure that in the sum, there are only finite terms. It follows that $\sum_{x \leq t \leq y} \mu(x,t) = 0$ for every $x < y$. From this definition we have at once $\mu(x,y) = 0$ if $x \nleq y$. For the order by divisibility we have $\mu(a,b) = \mu(b/a)$ where the second $\mu$ is the classical Moebius function. For the order by inclusion we have $\mu(T,U) = (-1)^{|U \setminus T|}$ for every $T \subseteq U$. For the order by exclusion we have $\mu(U,T) = (-1)^{|U \setminus T|}$ for every $U \supseteq T$. In this book we consider complex valued functions $\alpha(x,y)$ of two variables $x,y \in P$ as a matrices $|P| \times |P|$. Here $|P|$ denotes the cardinality of $P$. Moreover, we always assume that $\alpha(x,y) = 0$ if $x \nleq y$. These matrices (functions) are upper triangular. The product $\alpha \circ \beta$ is defined by $(\alpha \circ \beta)(x,y) = \sum_{x \leq t \leq y} \alpha(x,t)\beta(t,y)$. For example, the function $\delta(x,y) = \begin{cases} 1 & \text{if } x = y \\ 0 & \text{if } x \neq y \end{cases}$ is considered as identity matrix. A function $\alpha$ is called the inverse of $\beta$ if $\alpha \circ \beta = \beta \circ \alpha = \delta$. Now define $\xi(x,y) = \begin{cases} 1 & \text{if } x \leq y \\ 0 & \text{if } x \nleq y \end{cases}$ then $(\mu \circ \xi)(x,y) = \sum_{x \leq t \leq y} \mu(x,t)\xi(t,y) = \begin{cases} 1 & \text{if } x = y \\ 0 & \text{if } x \neq y \end{cases}$. Hence $\xi$ is inverse of $\mu$. More exactly, any finite poset $S$ can be numbered as $x_1, x_2, \cdots, x_n$ such that if $x_i < x_j$ then $i < j$. Thus, if $g(x) = \sum_{y \leq x} f(y)$ then

$$\begin{bmatrix} g(x_1) \\ g(x_2) \\ \vdots \\ g(x_n) \end{bmatrix} = \begin{bmatrix} 1 & 0 & \cdots & 0 \\ * & 1 & \cdots & 0 \\ \vdots & \vdots & \ddots & \vdots \\ * & * & \cdots & 1 \end{bmatrix} \begin{bmatrix} f(x_1) \\ f(x_2) \\ \vdots \\ f(x_n) \end{bmatrix} = (I+N) \begin{bmatrix} f(x_1) \\ f(x_2) \\ \vdots \\ f(x_n) \end{bmatrix}$$

where $N$ is a nilpotent matrix and $*$ may be 0 or 1. Clearly, $(I+N)^{-1} = I - N + N^2 - \ldots + (-1)^{n-1} N^{n-1}$ and we have

$$\begin{bmatrix} f(x_1) \\ f(x_2) \\ \vdots \\ f(x_n) \end{bmatrix} = \left( I - N + N^2 - \ldots + (-1)^{n-1} N^{n-1} \right) \begin{bmatrix} g(x_1) \\ g(x_2) \\ \vdots \\ g(x_n) \end{bmatrix}.$$

**Moebius inversion theorem**: If $f$ is a function from a finite poset $S$ into any commutative ring then $g(x) = \sum_{y \leq x} f(y)$ implies $f(x) = \sum_{y \leq x} g(y) \mu(y, x)$.

Now we choose $S = \{1, 2, \cdots, n\}$ and let $A_1, A_2, \cdots, A_n$ be finite sets. For the order by divisibility and get back classical inversion formulas. For the order by inclusion, we let $f(S) = 0$ and $f(T) = \left| \bigcap_{i \in T'} A_i \setminus \bigcup_{j \in T} A_j \right|$ for a proper subset $T$ of $S$. Then $g(T) = \sum_{U \subseteq T} f(U) = \left| \bigcup_{j \in T'} A_j \right|$ and $g(S) = \left| \bigcup_{j \in S} A_j \right|$ and we have the **inclusion–exclusion principle** $f(S) = \sum_{T \subseteq S} (-1)^{|T'|} g(T)$ or more exactly,

$$\left| \bigcup_{i=1}^{n} A_i \right| = \sum_{i=1}^{n} |A_i| - \sum_{i,j:1 \leq i < j \leq n} |A_i \cap A_j|$$
$$+ \sum_{i,j,k:1 \leq i < j < k \leq n} |A_i \cap A_j \cap A_k| - \cdots + (-1)^{n-1} |A_1 \cap \cdots \cap A_n|$$

But the induction according to $n$ is more natural proof. Let $A$ denote the union $\bigcup_{i=1}^{n} A_i$ of the sets $A_1, \ldots, A_n$. To prove the inclusion–exclusion principle in general, we first have to verify the identity

$$1_A = \sum_{k=1}^{n} (-1)^{k-1} \sum_{\substack{I \subseteq \{1,\ldots,n\} \\ |I| = k}} 1_{A_I} \quad (*)$$

for indicator functions, where

$$A_I = \bigcap_{i \in I} A_i.$$

There are at least two ways to do this:

**First possibility:** It suffices to do this for every $x$ in the union of $A_1, ..., A_n$. Suppose $x$ belongs to exactly $m$ sets with $1 \leq m \leq n$, for simplicity of notation say $A_1, ..., A_m$. Then the identity at $x$ reduces to

$$1 = \sum_{k=1}^{m}(-1)^{k-1} \sum_{\substack{I \subset \{1,...,m\} \\ |I|=k}} 1.$$

The number of subsets of cardinality $k$ of an $m$-element set is the combinatorical interpretation of the [binomial coefficient](#) $\binom{m}{k}$. Since $1 = \binom{m}{0}$, we have

$$\binom{m}{0} = \sum_{k=1}^{m}(-1)^{k-1}\binom{m}{k}.$$

Putting all terms to the left-hand side of the equation, we obtain the expansion for $(1 - 1)^m$ given by the [binomial theorem](#), hence we see that (*) is true for $x$.

**Second possibility:** The following function is identically zero

$$(1_A - 1_{A_1})(1_A - 1_{A_2}) \cdots (1_A - 1_{A_n}) = 0,$$

because: if $x$ is not in $A$, then all factors are $0 - 0 = 0$; and otherwise, if $x$ does belong to some $A_m$, then the corresponding $m$th factor is $1 - 1 = 0$. By expanding the product on the right-hand side, equation (*) follows. Now let $A_1, A_2, \cdots, A_n$ be events in the probability space. Integrate (*) according to the probability measure we have

$$\mathbb{P}\left(\bigcup_{i=1}^{n} A_i\right) = \sum_{k=1}^{n}(-1)^{k-1} \sum_{\substack{I \subset \{1,...,n\} \\ |I|=k}} \mathbb{P}(A_I),$$

or equivalently,

$$\mathbb{P}\left(\bigcup_{i=1}^{n} A_i\right) = \sum_{i=1}^{n}\mathbb{P}(A_i) - \sum_{i,j:i<j}\mathbb{P}(A_i \cap A_j)$$

$$+ \sum_{i,j,k:i<j<k}\mathbb{P}(A_i \cap A_j \cap A_k) - \cdots + (-1)^{n-1}\mathbb{P}\left(\bigcap_{i=1}^{n} A_i\right),$$

An **incidence algebra** is an [associative algebra](#), defined for any [locally finite poset](#) and [commutative ring](#) with unity. The members of the incidence algebra are the [functions](#) $f$ assigning to each nonempty interval $[a, b]$ a scalar $f(a, b)$, which is taken from the *ring of scalars*, a commutative [ring](#) with unity. On this underlying set one defines addition and

scalar multiplication pointwise, and "multiplication" in the incidence algebra is a
convolution defined by

$$(f * g)(a,b) = \sum_{a \leq x \leq b} f(a,x)g(x,b).$$

An incidence algebra is finite-dimensional if and only if the underlying poset is finite. An incidence algebra is analogous to a group algebra; indeed, both the group algebra and the incidence algebra are special cases of a categorical algebra, defined analogously; groups and posets being special kinds of categories. The multiplicative identity element of the incidence algebra is the **delta function**, defined by

$$\delta(a,b) = \begin{cases} 1, & \text{if } a = b \\ 0, & \text{if } a < b \end{cases}$$

# $A_n$, $n > 4$, Alternating groups

**Simplicity:** Solvable for $n < 5$, otherwise simple.

**Order:** $n!/2$ when $n > 1$.

**Schur multiplier:** 2 for $n = 5$ or $n > 7$, 6 for $n = 6$ or 7; see Covering groups of the alternating and symmetric groups

**Outer automorphism group:** In general 2. Exceptions: for $n = 1$, $n = 2$, it is trivial, and for $n = 6$, it has order 4 (elementary abelian).

**Other names:** $Alt_n$.

If $n \geq 3$ then $A_n$ is generated by the $3$-cycles
$(abc) = (bc)(ac) \ (ab)(cd) = (acd)(acb)$.

If $n \geq 5$ and $K$ is a normal subgroup of $A_n$ which contains a $3$-cycle then $K$ contains every $3$-cycle and consequently, $K = A_n$.

Now let $\alpha \in K \setminus \{\text{id}\}$ with maximal number of fix points. Then $\alpha$ is a $3$-cycle. If otherwise, $\alpha$ moves at least 4 numbers $1,2,3,4$, say. We write $\alpha$ as a product of disjoint cycles $\alpha = (123\cdots)\cdots$ or $\alpha = (12)(34)\cdots$. In the first case there is a cycle of length $\geq 3$ and in the second case every cycle is of length 2 (disjoint transpositions). Moreover, in the first case, $\alpha$ moves at least one other number 5, say (because $\alpha$ is even permutation). Now we have $\beta = (345) \in A_n$ and consequently, $\alpha_1 = \beta\alpha\beta^{-1} \in K$ ($K$ is a normal subgroup of $A_n$). Clearly, $\alpha_1 \neq \alpha$ so $\alpha_2 = \alpha_1\alpha^{-1} = \beta\alpha\beta^{-1}\alpha^{-1} \neq \text{id}$. Thus, $\alpha_2 \in K \setminus \{\text{id}\}$ and fixes 2. Moreover, if a number larger than 5 is fixed by $\alpha$ then it is also fixed by $\alpha_2$. Therefore, $\alpha_2$ fixes more numbers than $\alpha_2$ in the first case. In the second case $\alpha_2$ fixes 1. Thus, $\alpha_2$ fixes more numbers than $\alpha_2$ in any case. This is a contradiction. Therefore, $\alpha$ is a $3-\text{cycle}$ and consequently, $K = A_n$. Therefore, the symmetric group $S_n$ is not solvable if $n \geq 5$.

**Acknowledgement.** Deepest appreciation is extended towards the NAFOSTED (the **National Foundation for Science and Techology Development** in Vietnam) for the financial support.